\input amstex

%%%%%%%%%%%%%%%%%%%%%%%%%%%%%%%%%%%%%%%%%%%%%%%%%%%%%%%%%%%%%%%%%%%%%%%%%
%%%  Macros for this paper

\def\today{\ifcase\month\or
January\or February\or March\or April\or May\or June\or
July\or August\or September\or October\or November\or December\fi
\space{\number\day}, {\number\year}}

\def\s2{\sqrt 2}

\def\eps{\epsilon}
\def\sech{\mathop{\roman {sech}}\nolimits}

\def\Int#1{\raise 1pt\hbox{$\scriptstyle\int$}\!{#1}\,}
\def\Int#1{\kern 1pt\raise .15pt\hbox{$\char'177$}\!{#1}\kern .3pt}
\def\degree{\mathop{\roman {degree}}\nolimits}
\def\Der_#1(#2){#2\!_{{}_{(#1)}}}

\def\FT{\mathop{\Cal F}\nolimits}
\def\IFT
   {\mathop{\lower .25pt\hbox{$\Cal I$}\kern -2.1pt%
       \hbox{${\Cal F}$}}\nolimits}

\def\ICH
   {\mathop{{\Cal I}\kern -2.1pt\raise 0.15pt\hbox{${\Cal C\!H}$}}\nolimits}

\def\IST
   {\mathop{{\Cal I}\kern -2.1pt\raise 0.15pt\hbox{${\Cal S\!T}$}}\nolimits}

\def\Int{\mathop{\Cal J\kern -1pt}}

\def \Fr{\FT}
\def \IFr{\IFT}

\def\sph{\hbox{$\bold S$}}
\def\reals{\hbox{$\bold R$}}
\def\Cx{\hbox{$\bold C$}}

\def \norm|#1|{\left \Vert#1\right\Vert}
\def \Norm|#1|{\left \Vert \kern -0.75pt\left |#1\right\Vert\kern -0.75pt\right|}

\def\proof {\noindent{\it Proof.}\enspace}
\def\qedmark{\hbox{\vrule height 4pt width 3pt}}
\def\qedskip{\vrule height 4pt width 0pt depth 1pc}
\def\qed{\nobreak\quad\nobreak{\qedmark\qedskip}}
\def\HN{H_{\kern-1pt{}_N}}
\def\Hinfty{H_{\kern-1pt{}_\infty}}

\def\sphN{\sph_{\kern-1pt{}_N}}
\def\om{\omega}
\def\omN{\om_{\kern-1pt{}_N}}

\def\dsphN{{\hat{\sph}}_{\kern-1pt{}_N}}

\def \IP<#1,#2>{\left \langle #1,#2\right\rangle}
\def \PB(#1,#2){\left \lbrace #1,#2\right\rbrace}
\def \LD(#1?#2){{\Cal L}_{\!{}_{#1}}#2}

\def \L{\hbox{$\Cal L$}}

\def\a{{\bold a}}
\def\b{{\bold b}}
\def \mat(#1;#2;#3;#4)%
{\pmatrix{{#1}&{#2}\cr{#3}&{#4}}}

\def\dip<<#1;#2>>{\hbox{$<\!<$}#1,#2\hbox{$>\!>\!$}}

\def\diag(#1,#2){\hbox{{\rm diag}}(#1,#2)}

\def \mat(#1;#2;#3;#4)%
{\left(\matrix#1&#2\\#3&#4\endmatrix\right)}

\newcount\secnum \secnum=0       %%%% counter for section numbers %%%%
\newcount\subsecnum              %%%% counter for subsection numbers %%%%

\def\section#1{\advance\secnum by 1 \subsecnum=0% 
            \head{\the\secnum. #1}\endhead }

\def\subsection#1{\advance\subsecnum by 1%
       \subhead\nofrills{\the\subsecnum. #1.\ }\endsubhead}
\newcount\firstpageno
\newcount\tocpageno

\def \Dm{{\Cal D}_-}

\def \Dmd{\Dm^{\!\!^{\raise 0.75pt\hbox{$\scriptscriptstyle\roman {disc}$}}}}
\def \Dmc{\Dm^{\!\!^{\raise 0.75pt\hbox{$\scriptscriptstyle\roman {cont}$}}}}

%%%%  END OF MACROS FOR THIS PAPER
%%%%  BEGINNING OF DOCUMENT

\documentstyle{amsppt}

\def\qedmark{\hbox{\vrule height 4pt width 3pt}}
\def\qedskip{\vrule height 4pt width 0pt depth 1pc}
\def\qed{\nobreak\quad\nobreak{\qedmark\qedskip}}

%\NoBlackBoxes
\nologo

\hoffset = 52pt
\parskip = 0.3\baselineskip

\topmatter
\title{The initial value problem for weakly nonlinear PDE}\endtitle

\author Richard S. Palais\endauthor

\address%
Dept. of Mathematics,
University of California at Irvine
Irvine, CA 92697\endaddress

\email palais\@uci.edu\endemail


\abstract 
We will discuss an extension of the pseudospectral method developed 
by Wineberg, McGrath, Gabl, and Scott for the numerical integration 
of the KdV initial value problem. Our generalization of their algorithm 
can be used to solve initial value problems for a wide class of evolution 
equations that are ``weakly nonlinear'' in a sense that we will make precise. 
This class includes in particular the other classical soliton equations 
(SGE and NLS). As well as being very simple to implement, this method 
exhibits remarkable speed and stability, making it ideal for use with 
visualization tools where it makes it possible to experiment in real-time 
with soliton interactions and to see how a general solution decomposes into 
solitons. We will analyze the structure of the algorithm, discuss some 
of the reasons behind its robust numerical behavior, and finally describe 
a fixed point theorem we have found that proves that the pseudospectral 
stepping algorithm converges.
\endabstract

\thanks During the preparation of the initial version of this paper, the author 
was a member of The Institute for Advanced Study, in Princeton, New Jersey\endthanks

%\date July 1, 2000 \enddate

\date September 2014 \enddate
 

\endtopmatter

\document

\section{Introduction}

\noindent
  So-called pseudospectral methods are used for the numerical integration of 
evolution equations, using discrete Fourier transforms instead of finite differencing to 
evaluate spatial derivatives. An excellent early article is [FW]. A surprising 
fact is that these methods often work very well for nonlinear equations.  
The time-stepping for pseudospectral methods is accomplished by a classical 
differencing scheme that can in principle be either explicit or implicit, but 
for the usual stability reasons, an implicit method such as Crank-Nicolson
(the trapezoidal rule) is usually preferred. However, when the equation is
nonlinear, the solution of the implicit equations that arise can present a
problem. One approach is to employ split-stepping; use Crank-Nicolson plus
Gaussian elimination for the linear terms, but fall back to an explicit 
method for the nonlinear terms. An alternative approach, pioneered in [WMGS] 
and that we will refer to as the WMGS method, is to treat the linear and 
nonlinear terms together, write the implicit equation in fixed-point form, 
and then solve it by an iteration scheme. 

  WGMS originally developed their method to solve the initial value problems for the 
KdV and KP equations with periodic boundary conditions, and we became aware of their 
technique via an early version of  [LS], in which Yi Li and D.~H. Sattinger report on a 
modifiied WGMS algorithm. In this paper, we will discuss a generalization of the WGMS algorithm 
to treat  the initial value problem for a fairly broad class of evolutionary PDE that
are  ``weakly nonlinear'', in the sense that their nonlinear terms are a lower order
perturbation of the linear part (see below for a precise definition) and we will
prove a convergence theorem  for the iteration method  that is at the heart of the
WGMS algorithm.

\section {Weakly Nonlinear PDE of Evolution}

\noindent  
Let $U$ denote a finite dimensional complex inner product
space and $V$ a vector space of $U$-valued functions on $\reals$. Usually
we will work in a fixed orthonormal basis $(e_1,\ldots,e_n)$ for $U$ and
use it to identify $U$ with $\Cx^n$, so that elements $u$
of $V$ can be considered as $n$-tuples $(u^1,\ldots,u^n)$ of  
complex-valued functions. (If $n = 1$ we shall say we are in the 
{\it scalar case\/} and we then identify $u$ with $u^1$.) 

  We will specify $V$ more precisely later, but the elements of $V$ will admit 
derivatives up to a certain order, and they will in most cases be required 
to be $2 \pi$-periodic, in which case we shall also consider them as functions
on the unit circle in the complex plane. If $u(t)$ is a curve in $V$ we 
will also write $u(x,t)$ for $u(t)(x)$. As usual we think of $t$ as denoting 
time and $x$ as space. We denote by $D$ the differentiation operator
${\partial\over \partial x}$ and we also write $u^i_{x} = Du^i$, 
$u^i_{xx} = D^2 u^i$, etc., and of course $u_t ={\partial u\over \partial t}$.  

We will be considering 
``evolution equations'' of the form  $u_t = F(u)$, where $F:V \to V$ should be 
thought of as a vector field on $V$, and its form will be a smooth function 
(usually polynomial) of the $u^i$ and their derivatives, $Du^j,D^2 u^k,\ldots$. 
Usually $F(u)$ will the sum  of a ``dominant'' linear differential operator, 
and a nonlinear part that we can consider as a ``small perturbation'' of this 
linear part. By a linear differential operator on $V$ we will always mean 
an operator of the form $u \mapsto \L(u) = (\L^1(u),\ldots,\L^n(u))$ where 
$\L^i(u) = \sum_{j=1}^n  \L^i_j(D) u^j$. Here each $\L^i_j(X)$ is a polynomial
with constant coefficients in an indeterminate $X$. In the scalar case 
$\L(u) = \L(D)u$ and we will often use $\L(D)u$ to denote $\L(u)$ in the 
general case too. 

  The simplest kind of nonlinear operator that we shall consider is a zero order 
nonlinear operator, by which we will mean a map of the form 
$u \mapsto G(u) = (G^1(u),\ldots G^n(u))$, where $G^i(u)(x)= G^i(u^1(x),\ldots,u^n(x))$ 
and $G^i(Y_1, \ldots,Y_n)$ is either a constant coefficient polynomial on 
$\Cx^n$ or more generally an entire function of these variables 
(i.e., given by a power series that converges for all values of $(Y_1, \ldots, Y_n)$). 
Of course, care must be taken to make sure that if $u \in V$ then also $G(u) \in V$. 
When we come to the rigorous proofs, we will assume that $V$ is one of the 
Sobolev Hilbert spaces $H^m(\reals,U)$ for $m > {1\over 2}$, and since it is well-known 
that $H^m(\reals,\Cx)$ is a Banach algebras, it follows easily that $G$ is a smooth 
map of $H^m(\reals,U)$ to itself.  The most general kind of nonlinearity that we will 
consider will be one that can be  factored into a composition of the form $M(D) G(u)$ 
where $M(D)$ is a linear differential  operator as above and $G(u)$ is a zero order 
nonlinearity. 

  If $L(X) = \sum_{m=1}^\ell a_m X^m$ is a complex polynomial, then the differential 
operator $L(D)$ is called {\it formally skew-adjoint\/} if 
$\IP<L(D)u_1,u_2> = -\IP<u_1,L(D)u_2>$ whenever $u_1$ and $u_2$ are smooth maps of 
$\reals$ into $U$ with compact support. Here $\IP<u,v>$ denotes the $L^2$ inner product, 
i.e.,  $\IP<u,v> := \int_{-\infty}^\infty \IP<u(x),v(x)>\,dx$. Integration by parts 
shows that $D$ is skew-adjoint. Moreover an odd power of a formally skew-adjoint operator 
(and $i$ times an even power) is clearly again formally skew-adjoint, so  it follows that 
$L(D)$ is formally skew-adjoint if and only if the coefficients $a_m$ are real for $m$ odd 
and imaginary for $m$ even, i.e., if and only if $L(ik)$ is imaginary for all real $k$, 
and it is this last condition that we shall use.

\bigskip

\definition {Definition} A system of partial differential equation of the form:  
$$u^i_t = L^i(D) u^i +  M^i(D) G^i(u).\leqno{(WNWE)}$$
is called a {\it weakly nonlinear wave equation\/} if:
\item {1)} Each $L^i(D)$ is a formally skew-adjoint operator and the polynomials
           $L^i(X)$ all have the same degree, $\ell$,
\item {2)} $\degree M^i(X) < \ell$,
\item {3)} $G^i(0) = 0$, so that $u(x,t) \equiv 0$ is a solution of (WNWE).
\enddefinition

\noindent 
In what follows we will denote the minimum difference, 
$\ell - \degree M^i(X)$, by $q$.
For the most part, we will be dealing with the case $n = 1$, 
in which case we put  
$L = L^1$ and $M = M^1_1$, so $\ell = \degree L$ and 
$q =  \degree L -  \degree M$, and 
a weakly nonlinear wave equation has the form:
$$u_t = L(D) u + M(D) G(u).\leqno{(WNWE)}$$
\noindent
Two important examples are the Korteweg-deVries Equation: 
$$u_t = - u_{xxx} - u u_x = - D^3 u - {1\over 2} D(u^2),\leqno{(KdV)}$$ 
and the Nonlinear Schr\"odinger Equation: 
$$u_t = i u_{xx} + i |u|^2 u = i D^2 u + i |u|^2 u .\leqno{(NLS)}$$
In the former case we have $L(X) = -X^3, M(X) = -{1\over 2} X, G(X) = X^2$,
and in the latter, $L(X) = i X^2, M(X) = i, G(X) = |X|^2 X$.
 In the next section we will see that the Sine-Gordon Equation:

$$u_{tt } = u_{xx} + \sin u \leqno{(SGE)}$$ 
also can be regarded as a weakly nonlinear wave equation.

\section {The Sine-Gordon Equation}

\noindent
 A natural reduction of the linear wave equation to a system of first 
order PDE is ${\partial u \over \partial t} = {\partial v \over \partial x}$,
${\partial v \over \partial t} = {\partial u \over \partial x} $. 
in fact, ${\partial^2 \over \partial t^2} u = 
{\partial \over \partial t} {\partial \over \partial x} v =
{\partial \over \partial x} {\partial \over \partial t} v =
{\partial^2 \over \partial x^2} u$.

  This suggests that to find a representation of Sine-Gordon as a weakly 
nonlinear wave equation, we should start with systems of the form  
${\partial u \over \partial t} = {\partial v \over \partial x} + F(u,v)$,
${\partial v \over \partial t} = {\partial u \over \partial x} + G(u,v)$
or ${\partial  \over \partial t}(u,v) = L(u,v) + (F(u,v),G(u,v))$ where 
$F$ and $G$ are entire functions on $\Cx \times \Cx$, and of course $F(0,0) 
= G(0,0) = 0$. We will next show that with appropriate choice of $F$ and $G$ 
we do indeed get Sine-Gordon, and moreover that essentially the only other 
equations of the form $u_{tt } = u_{xx} + \Gamma(u)$ that  arise in this way 
are the Klein-Gordon equation, $u_{tt } = u_{xx} + u$, and the Sinh-Gordon
equation $u_{tt } = u_{xx} + \sinh u$. 

  Starting out as above,
${\partial^2 u\over \partial t^2} = 
{\partial \over \partial t} ({\partial v \over \partial x} + F(u,v) )=
{\partial \over \partial x} {\partial v \over \partial t} + F_1 {\partial u \over \partial t} 
+ F_2 {\partial v \over \partial t} =
{\partial^2 u\over \partial x^2} +
(F_1 F +F_2 G) + {\partial u \over \partial x} (G_1 + F_2) + 
{\partial v \over \partial x}(G_2 + F_1)$. For the latter to be of the form 
$u_{tt } = u_{xx} + \Gamma(u)$ we must have
${\partial  \over \partial v} (F_1F +F_2 G) = 0$, 
$G_1 = - F_2$, and $G_2 = -F_1$, in which case $u_{tt } = u_{xx} + \Gamma(u)$ with 
$\Gamma = F_1F + F_2G$.

  Next note that these conditions on $F$ and $G$ give $F_{11} = -G_{21} = -G_{12} = F_{22}$, 
or in other words, $F$ is a solution of the one-dimensional wave equation, and 
hence a sum of a left moving wave and a right-moving wave:
$F(u,v) = h(u+v) + k(u-v)$. Then using $G_1 = - F_2$, and $G_2 = -F_1$ it follows that
$G(u,v) = k(u-v) - h(u+v)$, where $h(0) = k(0) = 0$  in order to make  
$F(0,0) = G(0,0) = 0$. The condition ${\partial  \over \partial v} (F_1F +F_2 G) = 0$ now gives 
${\partial  \over \partial v}(h'(u+v)k(u-v) + h(u+v)k'(u-v)) = 0$ or 
$h''(u+v) k(u-v) = h(u+v)k''(u-v)$, 
or ${h''(u+v)\over h(u+v) } = {k''(u-v) \over k(u-v)}$. Since $u+v$ and $u-v$ are
coordinates, the only way the last relation can hold identically is for both sides to be a 
constant $\lambda$, i.e. $h'' = \lambda h$ and $k'' = \lambda k$. 

  If $\lambda$ is negative, say $\lambda = -\omega^2$, then since $h(0) = k(0) = 0$,
it follows that $h(u) = A \sin(\omega u)$ and $k(u) = B \sin(\omega u)$. 
If we choose $\omega = {1\over 2}$ and $A = B = 1$ we get 
$F(u,v) = \sin({u \over 2}+{v\over 2}) + \sin({u \over 2}-{v\over 2}) = 
2\sin{u \over 2}\cos{v\over 2}$ and similarly 
$G(u,v) = -2\cos{u \over 2}\sin{v\over 2}$, and this gives the 
system of partial differential equations 
${\partial u \over \partial t} = {\partial v \over \partial x} + 2\sin{u \over 2}\cos{v\over 2}$,
${\partial v \over \partial t} = {\partial u \over \partial x} -2\cos{u \over 2}\sin{v\over 2}$,
and we will leave it to the reader to check that if $(u,v)$ is a solution of this 
system, then $u$ is a solution of the Sine-Gordon equation. (Other choices of $A,B$, 
and $\omega$ lead to equations that can be transformed to the Sine-Gordon equation by 
a simple re-scaling of independent and dependent variables. Similarly taking $\lambda = 0$ 
gives the Klein-Gordon equation, and $\lambda$ positive gives Sinh-Gordon.)

 While this system of PDE for $u$ and $v$ is not in the form (WNWE), if we 
define $u^1 = u + v$ and $u^2 = u - v$, then $u^1$ and $u^2$ satisfy:
\smallskip
\centerline {$u^1_t = +u^1_x + 2 \sin({u^1 \over 2} - {u^2\over 2})$,}
\smallskip
\centerline {$u^2_t = -u^2_x + 2 \sin({u^1 \over 2} + {u^2\over 2})$.}
\smallskip
\noindent
which is manifestly in the form (WNWE), with $L^1(X) = D$, $L^2(X) = -D$, and $M^i(X) = 1$,
and  moreover we can recover $u$ from $u^1$ and $u^2$ by $u = {u^1 + u^2 \over 2}$.

  To simplify the exposition, we will from now on assume we are in the scalar case, $n = 1$ 
and that $G$ is a polynomial. The modifications needed for the general case are obvious.

\section {The Generalized WGMS Method (Heuristics)}

\noindent
Let us assume that for some particular example of (WNWE) we know that
there is a unique solution $u(t)$ with the initial condition $u(0) \in V$.
Let $\Delta t$ be close to zero, and let us look for a time-stepping algorithm
that, given a sufficiently good approximation to $u(t)$ as input will produce
an approximation to $u(t') = u(t + \Delta t)$ as output. If we integrate (WNWE)
with respect to $t$, from $t$ to $t'$, and use the trapezoidal rule to 
approximate the integrals on the right hand side, we find:
$$\eqalign{
u(t') - u(t) &= {\scriptstyle {\Delta t\over 2}} 
    L(D) [u(t) + u(t')]\cr
   &+ {\scriptstyle{\Delta t\over 2}}M(D)
       [G(u(t)) + G(u(t'))] \cr }$$
or 
$$\eqalign{
(I - d L(D)) u(t') &= 
(I  + d L(D)) u(t) \cr
   &+ d M(D)
       [G(u(t)) + G(u(t'))], \cr }$$
which we can rewrite as:
$$ u(t') = C u(t) + B[G(u(t)) + G(u(t'))]$$
where $d ={\Delta t\over 2}$, $B = {d M(D) \over I - d L(D)}$,
and  $C ={I + d L(D)\over I - d L(D)}$
is the Cayley transform of the skew-adjoint operator $d L(D)$.
We note that the skew-adjointness of $L(D)$ assures that
$I - d L(D)$ is invertible, and that $C$ is a unitary operator.
In fact, as we shall see shortly, on the Fourier transform side,
both $C$ and $B$ become simple multiplication operators, whose
properties are obvious from those of the polynomials $L(X)$ and $M(X)$.

  Next, for each $u$ in $V$, we define a map $H_u : V \to V$ by
$$H_u(w) := C u + B[G(u) + G(w)],$$ 
and we note that the equation above becomes $H_{u(t)}(u(t')) = u(t')$,
i.e., $u(t')$, which is what we are trying to compute, is a fixed-point
of $H_{u(t)}$.

  Now, we permit ourselves a little optimism---we assume that $u(t')$ is
in fact a {\it contracting\/} fixed point of $H_{u(t)}$. If this is so then,
for $\Delta t$ small, $u(t)$ will be close to $u(t')$, and we can expect
that iterating $H_{u(t)}$ starting at $u(t)$, will produce a sequence that 
converges to $u(t')$. This is the essence of the WGMS time-stepping algorithm 
(generalized to WNWE).

  For this to work as a numerical method, we must be able to compute
$H_u$ efficiently, and that is where the Fourier Transform comes in.
Let us write $\Fr$ for the Fourier Transform, mapping $V$ isomorphically 
onto $\hat V$, and $\IFr$ for its inverse. We define operators
$\hat C = \Fr C \IFr$ and $\hat B = \Fr B \IFr$ on $\hat V$. Then  
$\Fr H_u(w) = \Fr C \IFr \Fr (u) + \Fr B \IFr \Fr[G(u) + G(w)]$,
so we can rewrite $H_u$ as:
$$ H_u(w) = \IFr(\hat C \hat u + \hat B \Fr[G(u) + G(w)]),$$
where $\hat u = \Fr(u)$ is the Fourier Transform of $u$.

   Assuming that we have a good algorithm for computing $\Fr$ and $\IFr$
(e.g., the Fast Fourier Transform), it is now clear that it is easy and
efficient to calculate $H_u$, and hence to carry out the iteration. Indeed,
calculating $G(u)$ and $G(w)$ at a point $x$ is just a matter of evaluating 
the polynomial $G$ at $u(x)$ and $w(x)$. And since $M(X)$ and $L(X)$ are
constant coefficient polynomials, the operators $\hat C$ and $\hat B$ are
diagonal in the Fourier basis $e_k(x) = e^{ikx}$, i.e., they are multiplication
operators, by the rational functions ${1 + d L(i k)\over 1 - d L(i k)}$
and ${d M(i k) \over 1 - d L(i k)}$ respectively.  Since $L(D)$ is by
assumption skew-adjoint, $L(i k)$ is pure imaginary, so the denominator 
$1 - d L(i k)$ does not vanish. Moreover the function  
${1 + d L(i k)\over 1 - d L(i k)}$ clearly takes it values on the unit 
circle, and since $L(X)$ has degree greater than $M(X)$, it follows
that while the nonlinearity $G(w)$ may push energy into the high 
frequency modes of the Fourier Transform, multiplication by 
${d M(i k) \over 1 - d L(i k)}$ acts as a low-pass filter, attenuating
these high frequency modes and giving the WGMS method excellent numerical
stability.

\section {Proof that $H_u$ is a Contraction}

\noindent
In this section we will justify the above optimism by showing that, with 
a proper choice of the space $V$, a suitable restriction of the mapping $H_u$ 
does indeed satisfy the hypotheses of the Banach Contraction Theorem provided 
$\norm |u|$ and $\Delta t$ are sufficiently small. The space we will choose for 
$V$ is the Sobolev Hilbert space $H^m = H^m(\sph^1,V)$, with $m > {1\over 2}$.  
We recall that this is the Hilbert space of all functions $u$ in $L^2(\sph^1,V)$
such that $\norm |u|_m^2 = \sum_k (1 + k^2)^{m\over2}|\hat u(k)|^2$ is finite,
where as before, $\hat u(k)$ are the Fourier coefficients of $u$. 

  The principal property of these spaces that we shall need is that $H^m(\sph^1,\reals)$ 
is a commutative Banach algebra under pointwise multiplication when $m > {1\over 2}$ 
(cf. [A], Theorem 5.23, or [P]). As a first consequence, it follows that if
$P:V \to V$ is a polynomial mapping, then $u \mapsto P(u)$ is a map of 
$H^m$ to itself, and moreover $\norm |P(u)|_m < C \norm |u|_m^r$, where
$r$ is the degree of $P$. We will permit ourselves the abuse of notation of
denoting this latter map by $P$, and it is now elementary to see that
it is Frechet differentiable, and in fact that $DP_u(v) = P^\prime(u) v$, 
where $P^\prime$ is the derivative of $P$. (This will follow if we can show
that there is an algebraic identity of the form
$P(X + Y) = P(X) + P^\prime(X) Y + Q(X,Y) Y^2$, for some polynomial $Q$ in 
$X$ and $Y$. But it is clearly enough to check this for monomial $P$, 
in which case it is immediate from the binomial theorem.)  

  Let us denote by $B_R$ the ball of radius $R$ in $H^m$. Then as an immediate
consequence of the preceeding remarks we have:

\proclaim {Proposition 1} For any $R > 0$ there exist positive constants 
$C_1$ and $C_2$ such that $\norm |G(u)|_m < C_1$ and $\norm |DG_u| < C_2$ 
for all $u$ in $B_R$.
\endproclaim

It will be important for us to have a good estimate of how the norm of $B$ 
depends on $\Delta t$.

\proclaim {Proposition 2} Given $T > 0$, there is a positive constant $C_3$ 
such that the norm of the operator $B$ on $H^m$ satisfies 
$\norm |B| < C_3 {\Delta t}^{q\over\ell}$, for all $\Delta t < T$, where 
$\ell = \degree(L(X))$ and $q = \degree(L(X)) - \degree(M(X))$. Thus 
$\lim_{\Delta t \to 0} \norm |B| = 0$.
\endproclaim
\proof
It is clear that the Fourier basis $e_k(x) = e^{ikx}$ is orthogonal with respect 
to the $H^m$ inner-product (though not orthonormal, except for the case 
$H^0 = L^2$). Thus, since all constant coefficient differential operators are 
diagonalized in this basis, we can compute their norms on $H^m$ by taking the maximum
absolute values of their eigenvalues on the $e_k$. In the case of $B$, we have 
already seen that these eigenvalues are  ${d M(i k) \over 1 - d L(i k)}$.
Since $d = {\Delta T\over 2}$, to prove the proposition it will suffice to show
that ${d M(i k) \over 1 - d L(i k)} < C_3 d^{q\over \ell}$ for all real $k$  
and all $d < 2T$.

  Writing $L(X) = \sum_{j=0}^\ell b_j X^j$ and $M(X) = \sum_{j=0}^{\ell-q}a_j X^j$,
let us define parametric families of polynomials $L_c$ and $M_c$ for 
$c \ge 0$ by $L_c(X) = \sum_{j=0}^\ell (c^{\ell-j} b_j) X^j$ and 
 $M_c(X) = \sum_{j=0}^{\ell-q}(c^{\ell-q-j} a_j) X^j$.
Now note that if we define $\delta = d^{1\over \ell}$ then 
(since $\delta^{\ell-j} (\delta X)^j = d X^j$) 
clearly $L_\delta(\delta X) = d L(X)$, and similarly
$M_\delta(\delta X) = \delta^q d M(X)$, so 
${d M(i k) \over 1 - d L(i k)} = 
   \delta^q{M_\delta(i\delta k) \over 1 - L_\delta(i\delta k)},$
and to complete the proof it will suffice to show that the family
of rational functions 
$R_c(x) = {M_c(i x) \over 1 - L_c(i x)}$
is uniformly bounded for $0 \le c \le (\Delta T/2)^{1\over \ell}$ and 
$x$ real. If $\tilde\bold R$ is the one-point compactification of 
$\reals$ and we define $R_c(\infty) = 0$, then since the denominator
of $R_c(X)$ never vanishes and has degree greater than the numerator,
if follows that $(c,x) \mapsto R(c,x)$ is continuous and hence bounded on the 
compact space $[0, (\Delta T/2)^{1\over \ell}] \times \tilde \bold R$.
\qed 

\proclaim {Theorem} Given $R>0$ there exist positive $r$ and $T$ such that
$H_u$ is a contraction mapping of $B_R$ into itself provided that $u$ is 
in $B_r$ and $\Delta t < T$. Moreover there is a uniform contraction constant 
$K < 1$ for all such $u$ and $\Delta t$.
\endproclaim
\proof
From Proposition 1 and the definition of $H_u$ it follows that $H_u$ is 
differentiable on $H^m$ and that $D(H_u)_v = B\circ DG_v$. 
Then, again by Proposition 1, 
$\norm |D(H_u)_v| < C_2 \norm |B|$ for all $u$ in $B_R$,
and so by Proposition 2, 
$\norm |D(H_u)_v| < C_2 C_3 {\Delta t}^{q\over\ell}$.
Given $K < 1$, if we choose 
$T < {\left({K\over C_2 C_3}\right)}^{\ell\over q}$ then
$\norm |D(H_u)_v| < K$ on the convex set $B_R$ and hence $K$ is a 
contraction constant for $H_u$ on $B_R$, and it remains only to show 
that if we choose $r$ sufficiently small, and perhaps a smaller $T$ 
then $H_u$ also maps $B_R$ into itself for $u$ in $B_r$. 

   But using the definition of $H_u$ again, it follows that 
$$\norm |H_u(w)|_m < \norm |Cu|_m + 
  \norm |B|(\norm |G(u)|_m + \norm |G(w)|_m),$$
and recalling that $C$ is unitary on $H^m$, it follows from Propositions 1
and 2 that  $\norm |H_u(w)|_m < r + 2 C_1 C_3 T^{q\over \ell}$.
Thus $H_u$ will map $B_R$ into itself provided 
$r + 2 C_1 C_3 T^{q\over \ell} < R$, i.e., provided  
$r < R$ and $T < {\left({R-r\over 2 C_2 C_3}\right)}^{\ell\over q}$. 
\qed

  This completes a constructive proof of short-time existence for 
equations of the WNWE type. We note that for the standard examples, 
KdV and KP, it is well-known that the solutions exist for all time, and 
it would be interesting to know if this is true in general, and if not to
have a specific counter-example. 

\bigskip
 
\section {Numerics}

\noindent
There are several types of numerical errors inherent in the
WGMS algorithm. The first and most obvious is the error in 
approximating the integral of the right hand side of the equation
using the trapezoidal rule. A second ``truncation'' error occurs when 
we stop the fixed point iteration after a finite number of steps.
 
   In actually implementing the WGMS algorithm to solve an initial
value program numerically, one usually chooses an integer $N$ of 
the form $2^e$, and works in the space $V_N$ of ``band-limited'' 
functions $u$ whose Fourier coefficients $\hat u(k)$ vanish 
for $|k| > N/2$. Of course, $V_N$ is in all the Sobolev spaces. 
If we start with an initial condition $u_0$ not actually in 
$V_N$ then there will be an aliasing error when the initial 
Fast Fourier Transform projects it into $V_N$. Also, since the 
WGMS method does not rigorously preserve $V_N$, there will be 
further such errors at each time step.   It would be interesting to analyze 
these local errors and how they propagate in order to obtain a 
bound for the global error.

\bigskip\bigskip

\def\Bibliography{%
\def \bookitem //##1//##2//##3//##4//##5//##6//##7//##8//
%#1=key,#2=author,#3=title,#4=publisher,#6=year,#8=comment
{\ref%
\key \ignorespaces##1%
\by \ignorespaces##2%
\book \ignorespaces##3%
\publ \ignorespaces##4%
\yr \ignorespaces##6%
\finalinfo \ignorespaces##8%
\endref}
\def \b{\bookitem}
\def \articleitem //##1//##2//##3//##4//##5//##6//##7//##8//
%#1=key,#2=author,#3=title,#4=journal,#5=volume,#6=year,#7=pages,#8=comment
{\ref%
\key \ignorespaces##1%
\by \ignorespaces##2%
\paper \ignorespaces##3%
\jour \ignorespaces##4%
\vol \ignorespaces##5%
\yr \ignorespaces##6%
\pages \ignorespaces##7%
\finalinfo \ignorespaces##8%
\endref}
\def \a{\articleitem}
\def \preprintitem //##1//##2//##3//##4//##5//##6//##7//##8//
%#1=key,#2=author,#3=title
{\ref%
\key \ignorespaces##1%
\by \ignorespaces##2%
\paper \ignorespaces##3%
\toappear%
\endref}
\def \p{\preprintitem}
}
% end of bibliographic macros

\Refs

\widestnumber\key{AKNS2}

\Bibliography

\b //A//Adams,R.A.//Sobolev Spaces//Academic Press////1975//////

\a //FW//Forneberg,B. and Whitham, G.B.//A Numerical and Theoretical Study of 
Certain Nonlinear Wave Phenomena//Proc. R. Soc. Lond. A//289//1978//373--403////

\a //LS//Li, Y. and Sattinger, D.H.//Soliton Collisions in the Ion Acoustic Plasma Equations//\hskip 10pt J. math. fluid mech. //1//1999//117--130////

\b //P//Palais,R.S.//Foundations of Global Non-linear Analysis//
Benjamin and Co.////1968//////

\a //WGSS//Wineberg, S.B., Gabl,E.F., Scott, L.R. and Southwell, C.E.//Implicit 
Spectral Methods for Wave Propogation Problems//J. Comp. Physics//97//1991//311--336////

\endRefs

\section {Appendix---A MatLab Implementation}

  The following is a MATLAB M-file that implements the WGMS
algorithm to solve both KdV and NLS. It should be clear how to 
add further case statements to the definitions of the functions
Bhatfcn, G, and mfcn to extend this to handle other weakly nonlinear
equations.

%  VERBATIM MACROS  TeXbook  p380 et seq.
{\obeyspaces\global\let =\ } %let active space = control space
{\catcode`\`=\active\gdef`{\relax\lq}}

%  If <char> is any character and <text> 
%  is a string of characters that does not contain <char>, 
%  then \verbatim <char><text><char> 
%  defined below will set <text> verbatim 

%Here is a macro to allow you to input a complete file verbatim.
%Use \verbatiminput{filename} in much the same way as you would use
%\input filename.
% verbatim input macro
\def^^L{\par}
  
%  END VERBATIM MACROS

\verbatim@
function xxx = kdv_nls;
%KdV_nls.m:Korteveg de Vries & Nonlinear Schrodinger IVP
%Pseudo-spectral integration of
%u_t + uu_x + u_xxx = 0 and 
%iu_t+ u_xx + nu |u|^2 u =0
clear;
global KDVeqn NLSeqn TheEqn 
global nu disp a b m k dt;
global Cu B_hat G_u;
KDVeqn = 1; NLSeqn = 2;
TheEqn = input('Enter either 1 for kdv or 2 for nls ');
%initialization of variables for the two cases.
  switch TheEqn;
case NLSeqn
  N=1024;  period=2*pi;  nu = 2;  b = 3; 
  Top = b*sqrt(2/nu)*1.1; Bottom = - Top; 
  Left = - period/2; Right = period/2;
  TheEquationName = 'Non-Linear Schrodinger';
case KDVeqn
  N=512; period = 20;  disp = 0.05;
  Top = 1.75; Bottom = -0.1; 
  Left = - period/2; Right = period/2;
  TheEquationName = 'Korteweg de Vries';
end;
h = 2*pi/N;   				    % spatial increment (before scaling).
x = (-pi:h:pi - h);			    % unscaled space lattice
k = -i*[(0:N/2) (1-N/2:-1)];  % Fourier transform of d/dx.
a=period/(2*pi);			    % spatial scale-factor
t=0;						    % initial time
dt=0.01;					    % time step
y=a*x; 						    % scaled space lattice
u = initial_condition(x);		% initial condition
plothandle = plot(y,a*real(u));
set(plothandle,'erasemode','background');
axis([ Left Right Bottom Top]);
title(TheEquationName);
m = mfcn;
B_hat=bhatfcn;					 % Linear factor of nonlinear part of RHS
C_hat=(1+m)./(1-m);				 % Cayley transform of linear part of RHS
while 1;
Cu =  C_hat.*fft(u);
G_u = G(u);      				 %The nonlinearity of the RHS applied to u
w = u;
for n=1:3;
w = Iterator(w); % =ifft(Cu + B_hat.*fft(G_u + G(w)));
end
 u=w;
 set(plothandle,'ydata',a*real(u)); drawnow;  
 t=t+dt;
end

%Subsidiary functions used in this program

function nlnr = G(q);
global KDVeqn NLSeqn TheEqn;
  switch TheEqn
case NLSeqn
    nlnr = q.*abs(q).^2;
case KDVeqn
    nlnr = q.^2;
end;

function ic = initial_condition(x);
global KDVeqn NLSeqn TheEqn nu a b;
  switch TheEqn
case NLSeqn
    ic =  sqrt(2/nu)*b*sech(b*x);
case KDVeqn
    ic = exp(-1.2*a*x.^2)/a;
end;

function bhat = bhatfcn;
global KDVeqn NLSeqn TheEqn nu disp a b m k dt;
  switch TheEqn
case NLSeqn
    bhat = 0.5*i*nu*dt./(1-m);
case KDVeqn
    bhat = 0.5*dt*k./(1 - m); 
end;

function mfn = mfcn;
global KDVeqn NLSeqn TheEqn nu disp a k dt;
  switch TheEqn 
case NLSeqn

    mfn = 0.5*i*dt*k.^2;
case KDVeqn
    mfn = disp*a^(-3)*0.5*dt*k.^3; 
end;

function iter = Iterator(w);
global Cu B_hat G_u TheEqn;
  iter = ifft(Cu + B_hat.*fft(G_u + G(w)));
@

\enddocument

\bye